\documentclass{amsart}
\usepackage{amsmath, amssymb, graphics}
\usepackage{enumerate, amstext, amsthm}

\newtheorem{theorem}{Theorem}
\newtheorem{proposition}[theorem]{Proposition}

\newtheorem{example}[theorem]{Example}

\def\qed{\hfill\vbox{\hrule width 4 pt
\hbox{\vrule height 6 pt width 4 pt}} \medskip}
\def\IC{{\mathbb C}}

\def\diag{{\rm diag}\,}

\begin{document}
\title[Eigenvalue continuity and Ger\v{s}gorin's   theorem ]
{Eigenvalue continuity and Ger\v{s}gorin's  theorem}

\author{Chi-Kwong Li, Fuzhen Zhang}
\address[Li] {Department of Mathematics, College of William
\& Mary, Williamsburg, VA 13185, USA.}
\email{ckli@math.wm.edu}
\address[Zhang]{Department of mathematics, Nova Southeastern University,
Ft. Lauderdale, FL 33314, USA.}
\email{zhang@nova.edu}
\date{\today}


\maketitle

\begin{abstract}
Two types of eigenvalue continuity  are commonly used
in the literature. However, their meanings and the conditions under which continuities  are used are not always stated clearly. This can lead
 to some confusion and needs to be addressed. In this note,
we revisit  the
Ger\v{s}gorin disk theorem  and clarify the issue concerning the proofs of the theorem  by continuity.
\end{abstract}

\section{Introduction}

In his seminal paper in 1931 \cite{Ger31}, Ger\v{s}gorin presented an important result about the localization of the eigenvalues of matrices.  He showed that
(1) all eigenvalues of a square matrix lie in the union of the later so-called {\em Ger\v{s}gorin disks} and (2) if some, say $m$, of the disks are disjoint from the remaining disks, then the union of these $m$ disks contains exactly $m$ eigenvalues (counted with algebraic multiplicities). The result was named after  Ger\v{s}gorin as the {\em Ger\v{s}gorin disk theorem} due to its importance and applications for estimating  and  localizing   eigenvalues. 

Let $A=(a_{ij})$ be an $n\times n$ complex matrix and let $r_i=\sum_{j\not = i} |a_{ij}|$, $i=1, \dots, n$.
The set $D_i=\{z\in \Bbb C : |z-a_{ii}|\leq r_i\}$ is referred to as a Ger\v{s}gorin disk of $A$.
Let $A_0$ be the diagonal matrix that has the same main diagonal as  $A$.
Ger\v{s}gorin  proved the second  part of his theorem by considering the matrix $A(t)=A_0+t(A-A_0)$, $t\in [0, 1]$, and letting $t$ increase continuously  from 0 to 1. Intuitively,
the concentric Ger\v{s}gorin disks of $A(t)$ centered at
$a_{ii}$ ($i=1, \dots, n$)  get larger and larger  as $t$ increases from 0 to 1.
 He stated that
``Since the eigenvalues of the matrix depend continuously on its elements, it follows that $m$ eigenvalues must always lie in the disks ...".
Ger\v{s}gorin used   as a fact without justification that
{\em eigenvalues are continuous functions of the entries of matrices}.

Such a statement is often seen in the literature when it comes to the proof of the second part of the
Ger\v{s}gorin disk theorem.
For instance, here are a few widely-cited 
and comprehensive  references.   In the first edition of Horn and Johnson's book
{\em Matrix Analysis} \cite{HJ185}, page~345, it asserts that ``the eigenvalues are continuous functions of the entries of $A$ (see Appendix D)...", in Rahman and  Schmeisser's
{\em Analytic Theory of Polynomials} \cite{Rahman02}, page~55,
it states that
\lq\lq The eigenvalues of $A(t)$ are continuous functions of
$t$ ...\rq\rq,
in Varga's  {\em Ger\v{s}gorin and His Circles} \cite{Var2004}, page~8,
it is written that
``the eigenvalues $\lambda_i(t)$ of $A(t)$ also vary continuously
with $t$
...",
and
in Wilkinson's   {\em The Algebraic Eigenvalue Problem} \cite{Wil1965},
page~72,
 it says that
``the eigenvalues all traverse continuous paths".

What does it really mean to say that eigenvalues  are continuous functions?
 Ger\v{s}gorin's  proof by continuity may lead one to imagine {continuous}
curves of the eigenvalues  evolving 
on the complex plane,  or one may trace the curves continuously.
But that is not as easy as it sounds. First, ordering eigenvalues with a parameter
can be tricky and difficult; second,
the eigenvalue curves may merge 
 and the algebraic multiplicities of eigenvalues can change as the parameter varies.

\begin{example}\label{example1}
{\rm
Let $A(t) = \left( {0 \atop t}{t\atop 0}\right )$, $t\in [-1, 1].$  Each  $t$
produces a set of two eigenvalues. 
How does one order the eigenvalues as functions of $t$? It is natural to order the eigenvalues of $A(t)$ as $\lambda_1(t)=t$, $\lambda_2(t)=-t$.
Notice that $A(t)$ is real symmetric. For real symmetric (or complex Hermitian) matrices, we usually want  the eigenvalues to be in a non-increasing (or non-decreasing) order. So we would order the
eigenvalues as $\mu_1(t)=|t|\geq  \mu_2(t)=-|t|$. (Unlike $\lambda_1(t)$ and $\lambda_2(t)$, $\mu_1(t)$ and $\mu_2(t)$ are not differentiable.  Of course, there are
infinitely many  ways to parameterize the eigenvalues as non-continuous functions.)}
\end{example}\label{example0}

The eigenvalues in Example~\ref{example0}  are parameterized as continuous functions of $t$.
Is this always possible?
The answer is {\em yes} for $t$ on a real interval (but why? see Theorem \ref{thm1})
and {\em no} for $t$ on a complex domain containing the origin (see Example~\ref{example1}).

  Ger\v{s}gorin's original proof   by continuity is  more like  \lq\lq hand-waving" than a rigorous  proof and it
 has led to some confusion or ambiguity
\cite{Farenick2005}. A rigorous proof of the theorem using eigenvalues as continuous functions  requires creating or referencing some heavy machinery that was absent from all the classical sources. This issue deserves attention and clarification for both teaching and research.

Additionally, matrices depending on a parameter play   important roles in  scientific areas.
In some studies  such as stability problems and adiabatic quantum computing, one may consider a real parameter $t$  joining matrix $A$ and matrix  $B$ by $(1-t)A+tB$  and analyze the change of the eigenvalues
as $t$ varies.

In section 2, we briefly recap the eigenvalue continuity in the topological sense. In section 3, we summarize a celebrated result of Kato on the continuity of eigenvalues as functions. In section 4, we discuss the existing proofs of
the Ger\v{s}gorin disk theorem and present a proof with topological continuity and a proof with functional continuity. We  end the paper
by including  a short and neat proof of the second part of the Ger\v{s}gorin disk theorem by using   the argument principle.

\section{Topological continuity of eigenvalues}

Are eigenvalues of a matrix continuous functions of the matrix?
Since eigenvalue problems of matrices  are essentially  root problems of (characteristic) polynomials, one
immediately realizes that the question is a bit subtle and needs careful formulation.
It is  known that the roots of a polynomial vary continuously as {\em a function} of the coefficients. In \cite{Har_AMSP_1987}, the authors gave a nice proof for the
result concerning the continuity of  zeros of complex polynomials.
In fact, the map sending a monic polynomial
$f(z) = z^n + a_1z^{n-1} + \cdots + a_n$
to the multi-set of its zeros
$\pi(f) = \{\lambda_1, \dots, \lambda_n\}$
is continuous in the following   sense.
For   monic polynomials
$f(z) = z^n + a_1 z^{n-1} + \cdots + a_n$
and $\tilde f(z) = z^n + \tilde a_1 z^{n-1} + \cdots + \tilde a_n$
with multi-sets of zeros
$\pi(f) = \{\lambda_1, \dots, \lambda_n\}$ and
$\pi (\tilde f) = \{\tilde \lambda_1, \dots, \tilde \lambda_n\}$,
 one can  use  the metrics
$$\|f-\tilde f\| = \max\{|a_j-\tilde a_{j}|: 1 \le j \le n\}$$
and
{\small $$d(\pi (f), \pi (\tilde f)) =
\min_{J}
\left\{ \max_{1 \le j \le n}  |\lambda_j - \tilde \lambda_{j_\ell}|:
J=(j_1, \dots, j_n)
\hbox{ is a permutation of } (1, \dots, n)\right\}.$$}
Then $\pi$ is (pointwise) continuous; that is, for fixed $f$ and
for any given  $\varepsilon>0$, there exists  $\delta>0$ (depending on $f$)  such that
$\|f-\tilde f\| <\delta$ implies $d(\pi(f), \pi(\tilde f))<\varepsilon$.
Moreover, if $\xi$ is a zero of $f(z)$ with algebraic multiplicity $m$, then
$\tilde f$ has exactly $m$ zeros 
in the disk centered at $\xi$ with radius $\varepsilon$. 

If we identity $f$ with $n$-tuple $(a_1, \dots, a_{n})\in \IC^n$,
then $\pi$ is a homeomorphism between $\IC^n$ (with the usual topology) and the quotient space  $\IC^n_{\sim}$ (with the induced quotient topology),
the  unordered $n$-tuples (see \cite[p.\,153]{BhaMA97}) .

Applying this  result to matrices, one
gets the eigenvalue continuity  as
the eigenvalues of an $n\times n$ matrix $A$
 are the zeros of
the characteristic polynomial
$$p_A(z) = \det(zI - A) = z^n + a_1 z^{n-1} + \cdots + a_n,$$
where $a_j$ is $(-1)^j$ times the sum of the
$j\times j$ principal minors of $A$.

To be more specific, with $M_n$ for the space of $n\times n$ complex matrices,
we consider the eigenvalue function
$\sigma: M_n \rightarrow \IC^n_{\sim}$ that maps a matrix $A\in M_n$ to its spectrum
$\sigma (A)\in \IC^n_{\sim}$.
For the continuity of $\sigma$, we can use any (fixed) norm $\|\cdot\|$
on $M_n$.

The function $\sigma$ is continuous, i.e., for fixed $A\in M_n$ and for  any given $\varepsilon > 0$, there is
$\delta > 0$ (depending on $A$) such that
$d(\sigma(A), \sigma(\tilde A)) < \varepsilon$ whenever $\|A-\tilde A\| < \delta$.
Such an eigenvalue continuity
may be referred to as
{\em eigenvalue  topological  continuity} or {\em eigenvalue  matching continuity}.
Thus,  eigenvalues are always continuous in the topological sense.

  A nice proof regarding eigenvalue topological continuity for the discrete case (i.e., matrix sequences)  is available in \cite[pp.\,138--140]{Artin}.  The same continuity of eigenvalues is also studied
  in \cite[p.\,121]{HJ1_2nd}) by using Schur triangularization and compactness of the unitary group.
Closely related to eigenvalue continuity are
eigenvalue perturbation (variation) results with norm bounds
involving  the entries of matrices (see \cite{RCLiHandbookPert, Ost} and \cite[p.\,563, Appendix D]{HJ1_2nd}).

There is another possible way of thinking of the eigenvalue continuity problem. Let $A(t)$ be a family of $n\times n$ matrices depending continuously on a parameter $t$ over a domain in the complex plane or on a real interval. Then do there exist $n$ continuous complex functions of $t$  that   represent eigenvalues of $A(t)$?
We discuss the question in the next section.

\section{Parametrization of eigenvalues as continuous functions}

In some applications, one needs to consider a continuous
function $A: D \rightarrow M_n$, where $A(t) \in M_n$ and $D$ is a certain subset of
$\IC$ (say, a domain); and one wants to  parametrize the eigenvalues of $A(t)$ as
$n$ continuous functions $\lambda_1(t), \dots, \lambda_n(t)$ with $t \in D$.
We  refer to such continuity as {\em eigenvalue functional continuity} provided
there exist $n$ continuous functions of $t$  that represent eigenvalues of $A(t)$.

Eigenvalue functional continuity is widely used in
the proof of the second part of the Ger\v{s}gorin disk theorem;
similar ideas are needed in  the perturbation
theory of Hermitian matrices, stable matrices, etc.
However,
such a parametrization  is not always possible over a complex domain (\cite[p.\,154]{BhaMA97},
\cite[p.\,64;  p.\,108]{Kato95}).

\begin{example}\label{example1}
{\rm
 Let $A(t) = \left( {0 \atop t}{1\atop 0}\right )$, $t \in D = \{z\in \IC: |z| < 1\}.$
It is impossible   to have two continuous functions
$\lambda_1(t), \lambda_2(t)$ on $D$ representing the
eigenvalues of $A(t)$. This is because each eigenvalue
$\lambda$ of $A(t)$ satisfies $\lambda^2 = t$; thus, the desired
continuous functions $\lambda_1(t)$ and $\lambda_2(t)$ have
to satisfy $(\lambda_1(t))^2 = (\lambda_2(t))^2 = t$ for all $t$ on the open unit disk,
which is impossible (as is known, there is no continuous function $f$ on a disk  $D$ containing  the origin  such that $(f(z))^2=z$ for all $z\in D$).

However,
as $t\rightarrow 0$, $A(t)$ approaches $A(0)=\left( {0 \atop 0}{1\atop 0}\right )$
(entrywise),  which has repeated eigenvalue 0. Any small disk that contains the origin will contain two eigenvalues of  $A(t)$ when $t$ is close enough to $0$.
(This is what   topological continuity means.)}
\end{example}

 The difference between   topological continuity  and  functional continuity is that the eigenvalues (as a whole) are always topologically continuous but need not be continuous
  as individual functions.
 The two continuities for $A(t)$ are equivalent
  when the parameter  $t$  belongs to a real interval.
     This is well explained in \cite{BhaMA97, Kato95}.

In \cite[p.\,109, Theorem 5.2]{Kato95}, the following remarkable result is shown.

\begin{theorem}[Kato, 1966] \label{thm1} Suppose  that $D \subset \IC$ is a connected domain
and that $A: D \rightarrow M_n$ is a continuous function.
If   {\rm (1)}  $D$ is a  real interval, or  {\rm (2)}
$A(t)$   has only real eigenvalues, then there exist  $n$ eigenvalues $($counted with
algebraic multiplicities$)$ of $A(t)$ that
can be parameterized as continuous functions $\lambda_1(t),$ $ \dots,$ $ \lambda_n(t)$
from $D$ to $\Bbb C$.
In the second case, one can set
$\lambda_1(t) \ge \cdots \ge \lambda_n(t)$.
\end{theorem}

The study of eigenvalue functional continuity can be traced back at least as early as 1954
\cite{Rellich54}. 
Rellich \cite[p.\,39]{Rellich69} showed that individual eigenvalues are continuous functions when the matrices are Hermitian (in such case all eigenvalues are necessarily real). In his well-received book,
 Kato \cite[p.\,109]{Kato95}
 showed   that
 topological continuity   implies  functional continuity  when the parameter is restricted to a real interval or if all the eigenvalues of the matrices  are real, i.e., Theorem~\ref{thm1}.

It is tempting to extend Kato's result on a real interval for the parameter  to a domain (with
interior points) on the complex plane.  However, this is impossible.
Let $z_0\not = 0$ and let $D_{z_0}$ be an open disk centered at $z_0$ that does not contain the origin. Considering $A(z) = \left (
{0 \atop z-z_0} {\;1 \atop \;0} \right )$,
 $z\in D_{z_0}$,  we see that  there does not exist a continuous  eigenvalue function
of $A(z)$  on $D_{z_0}$.
Suppose, otherwise,  there is a  continuous eigenvalue function $\lambda(z)$ on  $D_{z_0}$, then
$(\lambda(z))^2 = z-z_0$ for all $z\in D_{z_0}$. This leads to a continuous function
$f(z) = \lambda(z+z_0)$ defined on the open unit disk $D=\{z\in \IC : |z|<1\}$  such that $(f(z))^2 = z$
for all $z\in D$, a contradiction.

So, in a sense,  the result of Kato is the best possible with respect to eigenvalue functional continuity.

\section{Proofs of the Ger\v{s}gorin disk theorem}

Ger\v{s}gorin's  disk theorem
is a useful result for estimating  and  localizing  the eigenvalues of a matrix. Usually and traditionally, the second part of the theorem is proved by considering the matrix $A(t)=A_0+t(A-A_0)$ (where $A_0$ is the diagonal matrix that has the same main diagonal as  $A$)  and by using eigenvalue continuity (see
\cite{Bru1994},  \cite[p.\,23]{HoodThesis17}, \cite[p.\,345]{HJ185}, \cite[p.\,74]{Hou1964}, \cite[p.\,372]{Lan1985},
\cite[p.\,499]{Meyer2000}, \cite[p.\,55]{Rahman02}, \cite[p.\,169]{SS1990},
   \cite[p.\,8]{Var2004},   \cite[p.\,72]{Wil1965}), and \cite[p.\,70]{ZFZbook11}).
   However, in these references, it is not always clear
    which types of  eigenvalue continuity conditions
    were used. If it is topological continuity,
    then one needs to add some details in the proofs to justify why the total number of eigenvalues in an isolated  region remains the same when $t$ increases from 0 to 1 (note that the algebraic multiplicity of an eigenvalue may change); if
    it is functional continuity (which is the case in most texts), then it would
    be nice to state Kato's result (or other references) as evidence of
    the existence of continuous functions that represent the eigenvalues.
 In the following, we
 state the Ger\v{s}gorin disk theorem and give two different proofs. One
 (Proposition~5) uses eigenvalue functional continuity (Kato's theorem) and exploits the fact that a continuous function takes a connected set into a connected set; the other (Proposition~6)
uses eigenvalue topological  continuity and exploits the fact that a continuous
function on a compact set is {\em uniformly continuous}.
In the latter, we completely avoid the continuity of each
eigenvalue as a function.

\begin{theorem}[Ger\v{s}gorin \cite{Ger31}, 1931]\label{GDT}
  Let $A = (a_{ij}) \in M_n$ and define the disks
$$D_i = \big \{ z \in \IC: |z-a_{ii}| \le \sum_{j\ne i} |a_{ij}|\big  \}, \quad i = 1, \dots, n.$$
Then \,
 {\rm (1)} All eigenvalues of $A$ are contained in the union $\cup_{i=1}^n D_i$.
{\rm (2)}  If $\cup_{i=1}^n D_i$ is the union of
$k$ disjoint connected regions $R_1, \dots, R_k$,
and $R_r$ is the union of $m_r$ of the disks $D_1, \dots, D_n$, then
$R_r$  contains exactly $m_r$ eigenvalues of $A$, $r = 1, \dots, k$.
\end{theorem}

Part (1) says that every eigenvalue of $A$ is contained in a Ger\v{s}gorin disk.
Its proof is easy,  standard, and  omitted here.
Part (2)  is immediate from Propositions \ref{Prop5} and \ref{Prop6}.
We call a union of some Ger\v{s}gorin disks a {\em Ger\v{s}gorin region}
(which in general need not be connected). In particular, the singletons of  diagonal entries are
degenerate  Ger\v{s}gorin regions. By a {\em curve} we mean the image (range) of
a continuous map from a real closed interval to the complex plane $\gamma:  [a, b] \mapsto \IC$.

\begin{proposition}\label{Prop5}
Let $A = (a_{ij}) \in M_n$ and  let $A(t)=A_0 + t(A-A_0)$,
where $t \in [0,1]$ and
$A_0 = \diag(a_{11}, \dots, a_{nn})$. Then each  continuous eigenvalue
curve of
$A(t)$ lies entirely in a connected Ger\v{s}gorin region of $A$.
\end{proposition}
\proof
By Kato's result (Theorem \ref{thm1}),  there exists a selection of
$n$  eigenvalues $\lambda_1(t), \dots, \lambda_n(t)$ of $A(t)$ that are continuous functions in $t$ on the real interval $[0, 1]$. Moreover,   part  (1) of the Ger\v{s}gorin disk theorem ensures that
  $\lambda_1(t), \dots, $ $\lambda_n(t)$ are contained
  in $\cup_{i=1}^k R_i$ for every $t\in [0, 1]$, and each set $\lambda_j([0, 1])$ is connected.

Let $r\in \{1, \dots, k\}$. Since  $R_r$ comprises $m_r$ disks (not necessarily different)  whose centers are $m_r$ elements of
the diagonal matrix $A_0$,  $m_r$ of the continuous eigenvalue curves $\lambda_1(t), \dots, \lambda_n(t)$ are  in $R_r$ at $t=0$. If $\lambda_j(0)\in R_r$, then
the connected set $$\lambda_j([0, 1])= \lambda_j([0, 1])\cap
\cup_{i=1}^k R_i = (\lambda_j([0, 1])\cap R_r) \cup(\lambda_j([0, 1])\cap    \cup_{i\not = r}^k R_i)$$
is the union of two disjoint closed sets, the first of which is nonempty. Therefore, the second set is empty and hence $\lambda_j([0, 1])\subset  R_r$.
\qed

The following proposition  considers the eigenvalues as a whole in a Ger\v{s}gorin region rather than focusing on an individual eigenvalue as a function. That is, we use eigenvalue topological continuity and avoid entirely (the difficult issue of) eigenvalue functional  continuity (which is not needed) to prove the assertion. 

\begin{proposition}\label{Prop6}
Let $A = (a_{ij}) \in M_n$ and let $A(t)=A_0 + t(A-A_0)$, where $t \in [0,1]$ and
$A_0 = \diag(a_{11}, \dots, a_{nn})$. Then  a connected Ger\v{s}gorin region of $A$ contains the same number of eigenvalues of $A(t)$ for all  $t\in [0, 1]$.
\end{proposition}

 \proof
 Every entry of $ A(t)$ is a continuous function of $t \in [0,1]$ and
  each Ger\v{s}gorin disk of $A(t)$  ($0\leq t \leq 1$) is contained in a Ger\v{s}gorin disk of $A=A(1)$ with the same corresponding center.  Let $R_r$, $r=1, \dots, k$,  be the connected Ger\v{s}gorin regions
 of $A$. (The number of connected Ger\v{s}gorin regions  for $A(t)$ may vary depending on $t$.)
Suppose that $R_r$ contains $m_r$ diagonal entries of $A$, i.e., $m_r$ eigenvalues of
$A_0=A(0)$ (counted with algebraic multiplicities). We claim that $R_r$ contains $m_r$ eigenvalues of $A(t)$ for all
$t\in [0, 1]$ (that is, the sum of the algebraic multiplicities of the eigenvalues
of $A(t)$ remains constant on each connected Ger\v{s}gorin region of $A$  as $t$ varies from 0 to 1).

Since the eigenvalues are topologically continuous over the compact set $[0, 1]$,
the continuity is uniform. To be precise,
 the map $\varphi: [0, 1]\mapsto \IC^n_{\sim}$ defined by $\varphi (t)=\sigma (A(t))$
 is uniformly continuous.

Let $\varepsilon > 0$ be such that $|x-y| > 2 \varepsilon$ for all
$x, y$ lying in any two disjoint Ger\v{s}gorin regions of $A$.  There is $\delta > 0$ (depending only on $\varepsilon $) such that
for any $t_1$ and $t_2$ satisfying  $0 \le t_1<t_2  \le 1$ and $t_2-t_1 < \delta$, the eigenvalues of
$A(t_1)$ and $A(t_2)$ can be labeled as
$\lambda_1,\dots, \lambda_n$ and $\mu_1, \dots, \mu_n$ such that
$|\lambda_j - \mu_j| < \varepsilon$ for  $j=1, \dots, n$.

We divide the interval $[0, 1]$ into $N$ subintervals:
 $0 = t_0 < t_1 < \cdots < t_N = 1$
such that $t_{i}-t_{i-1} < \delta$ for $i = 1, \dots, N$.
We  show   that  on each of the intervals
  $A(t)$ has $m_r$
eigenvalues in $R_r$ for $t\in [t_{i}, t_{i+1}]$, $i=0, 1, \dots, N-1$.

By assumption, $A(t_0)=A_0=A(0)$ has exactly $m_r$ eigenvalues in $R_r$.
For any $t\in [t_0, t_1]$, since   $t-t_0<\delta$,  $A(t)$   has exactly $m_r$ eigenvalues each of which is located in some disk centered at an eigenvalue of $A(t_0)$ with radius $\varepsilon$. By our choice of $\varepsilon$, all the $m_r$ eigenvalues of $A(t)$ are contained in
$R_r$ (i.e., not in other regions) for all $t\in [t_0, t_1]$.

Because $t_2-t_1<\delta$, for any  $t\in [t_1, t_2]$, $A(t)$ has exactly  $m_r$ eigenvalues that are close (with respect to $\varepsilon$)
to the $m_r$ eigenvalues of $A(t_1)$ in $R_r$.
Again, by our choice of  $\varepsilon$, all these
$m_r$ eigenvalues of $A(t)$ are also contained in $R_r$ 
for all $t\in [t_1, t_2]$.

Repeating the arguments for
$[t_2, t_3],$ $ \dots, $ $ [t_{N-1}, t_N]$,
we see that $A(t)$ has exactly  $m_r$ eigenvalues
in $R_r$ for each $t\in [0, 1]$. Thus,
$A(t_N)=A(1)=A$ has
exactly $m_r$ eigenvalues in the region $R_r$.   \qed

\section{A proof of Ger\v{s}gorin theorem  using the argument principle}\label{ComplexAnalysis}
The Ger\v{s}gorin  disk theorem is a statement about counting eigenvalues according to their algebraic multiplicities;
 it is essentially about counting zeros of a polynomial that depends on a parameter.
Thus,   Rouch\'{e}'s theorem would be a much  more natural and effective tool since it focuses squarely on what the theorem says about numbers of eigenvalues.
This approach does  not require the parameter $t$ to be real and it does not need the concept of any eigenvalue  continuity, functional or topological.

There is a short and neat proof of the second part of the Ger\v{s}gorin disk theorem that uses    the argument principle.  This approach was adopted in the second edition of Horn and Johnson's book {\em Matrix Analysis}
\cite[p.\,389]{HJ1_2nd}
 (see   also \cite[p.\,103]{Ser10}), while the proof by continuity used in the first edition \cite[p.\,345]{HJ185} was abandoned.

Let $\Gamma$ be a simple contour in the complex plane that surrounds the Ger\v{s}gorin region to be considered.
Let $p_{t}(z)$ be the characteristic polynomial of $A(t)$ for each given $t\in [0, 1]$.  By  the argument principle
\cite[p.\,123]{JBC78},
the number of zeros (counted with algebraic multiplicities)  of $p_{t}(z)$ inside $\Gamma$ is
$$m(t)=\frac{1}{2\pi i}
\oint_{\Gamma} \frac{p_{t}'(z)}{p_{t}(z)}dz.$$
On the other hand,
$f(t, z):=\frac{p_{t}'(z)}{p_{t}(z)}$
is a continuous function from $[0, 1]\times \Gamma $ to $\Bbb C$. By Leibniz's rule
\cite[p.\,68]{JBC78}, $m(t)$ is a continuous function on $[0, 1]$. As $m(t)$ is an integer, it has to be a constant. Thus,
$m(0)=m(1)$, which is the number of eigenvalues of $A$ in the Ger\v{s}gorin region.

Similar ideas
  using Rouch\'{e}'s theorem or winding numbers have been employed in the study of localization for
  nonlinear eigenvalue problems \cite{BinH13, BinH15, HoodThesis17}.
Eigenvalues as  functions deserve study and it is
an interesting (and classical) problem.  There is  a well-developed theory on the
smoothness of roots of polynomials
(see  \cite{AKLM98},
\cite[Chap.\,II, \S 4]{Kato95}, \cite{KLM04}, and
\cite{LR07}).
\bigskip

\noindent
{\bf Acknowledgment.}
This work was initiated by a talk given by the second author at the
International Conference on Matrix Theory with Applications
- Joint meeting of \lq\lq International Research Center for Tensor and Matrix Theory (IRCTM)", Shanghai University, China,
and  \lq\lq Applied Algebra and Optimization Research Center (AORC)",  Sungkyunkwan University, South Korea, at the
Shanghai University, December 17-20, 2018. Both authors thank the Centers for hospitality during the meeting and  thank
Roger Horn for  helpful discussions.


\begin{thebibliography}{99}

\bibitem{Artin} M. Artin.
\newblock
{\em Algebra}, 2nd edition.
\newblock Prentice Hall, 2011.


\bibitem{AKLM98}
D. Alekseevsky,  A. Kriegl,  P.W. Michor, and M. Losik.
\newblock  {Choosing roots of polynomials smoothly.}
\newblock {\em Israel J. Math.} 105 (1998) 203--233.

\bibitem{BhaMA97} R. Bhatia.
\newblock
     {\em Matrix Analysis.}
    \newblock
    Springer-Verlag, New York, 1997.

\bibitem{BinH13} D. Bindel and A. Hood.
  \newblock
Localization theorems for nonlinear eigenvalue problems.
  \newblock
{\em SIAM J. Matrix Anal. Appl.} 34 (2013) 1728--1749.

\bibitem{BinH15} D. Bindel and A. Hood.
  \newblock
Localization theorems for nonlinear eigenvalue problems.
  \newblock
{\em SIAM Review} 57 (2015) 585--607.

\bibitem{Bru1994} R. Brualdi and S. Mellendorf.
  \newblock
{Regions in the complex plane containing the eigenvalues of a matrix.}
  \newblock
{\em Amer. Math. Monthly} 101 (1994), no. 10, 975--985.

\bibitem{JBC78} 
J.B. Conway.
  \newblock
{\em  Functions of One Complex Variable} (Graduate Texts in Mathematics  11),
  \newblock  2nd Edition. Springer, 1978.

\bibitem{Farenick2005}
D. Farenick.
  \newblock {Book review: Ger\v{s}gorin and His Circles by R. Varga}.
  \newblock
{\em Notes of the Canadian Mathematical Society,} Vol. 37, No. 4, pp.~7--8, May 2005
(https://cms.math.ca/notes/v37/n4/Notesv37n4.pdf).

\bibitem{Ger31}
S. Gerschgorin.
  \newblock  {\"{U}ber die Abgrenzung der Eigenwerte einer Matrix}.
    \newblock
{\em Izv. Akad. Nauk SSSR Ser. Mat.} 1 (1931) 749-754.

\bibitem{Har_AMSP_1987}
G. Harris and C.  Martin.
  \newblock
{The roots of a polynomial vary continuously as a function of the coefficients.}
  \newblock
{\em Proc. Amer. Math. Soc.}
100 (1987), no. 2, 390--392.


\bibitem{HoodThesis17} A. Hood.
  \newblock
Localizing the eigenvalues of matrix-valued functions: analysis and applications.
  \newblock
Ph.D. dissertation, Cornell University, 2017.

 \bibitem{HJ1_2nd}   R.A. Horn and C.R. Johnson.
   \newblock
    {\em Matrix Analysis,} 2nd edition.
      \newblock
 Cambridge University Press,
     New York,  2013.

\bibitem{HJ185}   R.A. Horn and C.R. Johnson.
  \newblock
    {\em Matrix Analysis.}
      \newblock
 Cambridge University Press,
     New York, 1985.

\bibitem{Hou1964}
  A.S.  Householder.
    \newblock
  {\em  The Theory of Matrices in Numerical Analysis.}
    \newblock Reprint of 1964 edition. Dover Publications, New York, 1975.

\bibitem{Kato95} T. Kato.  \newblock
{\em Perturbation Theory for Linear Operators.}
  \newblock
Springer, 1995 (earlier editions in 1966, 1976).

   \bibitem{KLM04}
A. Kriegl,  M. Losik, and P.W. Michor.   \newblock
  {Choosing roots of polynomials smoothly. II.}
     \newblock {\em Israel J. Math.} 139 (2004) 183--188.

\bibitem{Lan1985} P. Lancaster and M. Tismenetsky.   \newblock
{\em The Theory of Matrices,} 2nd edition.
  \newblock  Computer Science and Applied Mathematics. Academic Press,   Orlando, 1985. 

\bibitem{RCLiHandbookPert}
R.-C. Li.   \newblock
{\em Matrix Perturbation Theory} in
 Handbook of Linear Algebra edited by L. Hogben,
 2nd edition.
  \newblock
CRC Press, \nolinebreak 2014.

\bibitem{LR07}
M. Losik and A.  Rainer.   \newblock
  {Choosing roots of polynomials with symmetries smoothly.}
    \newblock
{\em Rev. Mat. Complut.} 20 (2007), no. 2, 267--291.

\bibitem{Meyer2000} C. Meyer.   \newblock
{\em Matrix Analysis and Applied Linear Algebra}.
  \newblock  SIAM, Philadelphia, 2010.  


\bibitem{Ost} A.M. Ostrowski.  \newblock
{\em Solution of Equations and Systems of Equations},
  \newblock  2nd edition.
(Appendix K, Continuity of the fundamental roots as functions of  the elements of the matrix).
 \newblock
 Academic Press, New York and London, 1966.

\bibitem{Rahman02}
Q. Rahman and G. Schmeisser.    \newblock
{\em Analytic theory of polynomials.}
  \newblock
London Math Soc. Monographs. New Series, 26. The Clarendon Press, Oxford University Press, Oxford, 2002.

\bibitem{Rellich54}
F. Rellich.
  \newblock
{\em Perturbation Theory of Eigenvalue Problems.}
  \newblock
Courant Institute of Math Sciences, New York University, 1954.

\bibitem{Rellich69}
F. Rellich.  \newblock
{\em Perturbation Theory of Eigenvalue Problems}
(Assisted by J. Berkowitz).
  \newblock
    Gordon and Breach Science Publishers, New York-London-Paris, 1969. 

\bibitem{Ser10} D. Serre.  \newblock
 {\em Matrices: Theory and Applications}, 2nd edition.
   \newblock
Graduate Texts in Math 216,  Springer, New York, 2010.

\bibitem{SS1990} G.W. Stewart and J.-G. Sun. {\em Matrix Perturbation Theory.} Academic Press, Boston,   1990.  

\bibitem{Var2004}
R.S. Varga.   \newblock
 {\em Ger\v{s}gorin and His Circles.}
   \newblock Springer Series in Computational Mathematics, 36. Springer-Verlag, Berlin, 2004. 

\bibitem{Wil1965}
J.H. Wilkinson.  \newblock   {\em The Algebraic Eigenvalue Problem.}
  \newblock Monographs on Numerical Analysis. Oxford Science Publications. The Clarendon Press, Oxford University Press, New York, 1965;  Rev. ed. 1988.

\bibitem{ZFZbook11} F. Zhang.  \newblock
 {\em Matrix Theory: Basic Results and
Techniques}, 2nd edition.   \newblock
Springer, New York,  2011.

\end{thebibliography}
\end{document}